\documentstyle{amsppt}
\magnification 1200
\parskip=\medskipamount
\NoBlackBoxes 
\def\ve{\varepsilon} \def\vp{\varphi}
\def\arrowk{^\to{\kern -6pt\topsmash k}}
\def\arrowK{^{^\to}{\kern -9pt\topsmash K}}
\def\arrowr{^\to{\kern-6pt\topsmash r}}
\def\arrowvp{^\to{\kern -8pt\topsmash\vp}}
\def\arrowf{^{^\to}{\kern -8pt f}}
\def\arrowg{^{^\to}{\kern -8pt g}}
\def\arrowu{^{^\to}a{\kern-8pt u}}
\def\arrowt{^{^\to}{\kern -6pt t}}
\def\arrowe{^{^\to}{\kern -6pt e}}
\def\tk{\tilde{\kern 1 pt\topsmash k}}
\def\td{\ \bar{\kern-2pt \topsmash d}}
\def\barm{\bar{\kern-.2pt\bar m}}
\def\barN{\bar{\kern-1pt\bar N}}
\def\barA{\, \bar{\kern-3pt \bar A}}

\def\mod{\text {mod }}
\def\average{\lower 1pt\hbox{ $\sim $} \kern -9.5pt\topsmash \sum}
\hsize = 6.2true in

\vsize=8.2 true in \NoRunningHeads
\TagsOnRight
\topmatter
\title
{\smc Integral Apollonian circle packings and prime curvatures}
\endtitle
\author
J.~Bourgain\footnote{\hbox{The author was partially supported by NSF grants DMS-0808042 and DMS-0835373}}
\endauthor
\address
Institute for Advanced Study, Princeton, NJ 08540
\endaddress
\email
bourgain\@ias.edu
\endemail
\endtopmatter
\document

\noindent
{\bf Summary:} It is shown that any primitive integral Apollonian circle packing captures a fraction of the prime numbers.
Basically the method consists in applying the circle method, considering the curvatures produced by a well-chosen family of binary
quadratic forms.
\bigskip

\noindent
{\bf Introduction}

In this paper, we pursue a line of research initiated in [GLMWY] and [S] on the arithmetical properties of integral Apollonian circle packings (ACP for short)
in the plane.  
The reader is also referred to [B-F1] for certain background material.

Throughout the paper, we consider bounded ACP's which are primitive, meaning that all curvatures of the circles in the packing do not share a factor
greater than one.
Let us recall that the set of curvatures in a given packing $P$ is obtained by action of the Apollonian group $A$ on the root quadruple $(a, b, c, d)$ of
co-prime integers $a<0\leq b\leq c\leq d$, $a+b+c\geq d$.
The group $A$ is a subgroup of the orthogonal group associated to the Descartes quadratic form
$$
Q(x_1, x_2, x_3, x_4)= 2(x_1^2+x_2^2+x_3^2+x_4^2)-(x_1+ x_2+x_3+x_4)^2
$$
whose vanishing is tantamount with $x_1, x_2, x_3, x_4$ being curvatures of mutually tangent circles.
The group $A$ is generated by the matrices
$$
S_1=\pmatrix -1&2&2&2\\ 0&1&0&0\\ 0&0&1&0\\ 0&0&0&1\endpmatrix \qquad S_2= \pmatrix 1&0&0&0\\ 2&-1&2&2\\ 0&0&1&0\\ 0&0&0&1\endpmatrix
$$
$$
S_3= \pmatrix 1&0&0&0\\ 0&1&0&0\\ 2&2&-1&2\\
0&0&0&1\endpmatrix \qquad S_4 = \pmatrix 1&0&0&0\\ 0&1&0&0\\ 0&0&1&0\\ 2&2&2&-1\endpmatrix
$$

The basic problem is to describe the set of curvatures appearing in a given packing $P$; the ultimate  
hope is to establish a local to global principle cf. [F-S].
More modestly, [GLMWY] put forward the `positive density' conjecture, according to which the set of curvatures in an ACP form a subset of $\Bbb Z$ of positive
density.
Following up on a technique proposed by P.~Sarnak, this problem was solved affirmatively in [B-F1] (a slightly stronger result is obtained in this paper;
see Theorem 1 and the Remark following its proof).
Using a result due to Iwaniec on representing shifted primes by binary quadratic forms, Sarnak also pointed out that any primitive ACP produces at least
$ c\frac X{(\log X)^{3/2}}$ distinct prime curvatures at most $X$, for $X\to\infty$.
Based on new results on the representation of integers by binary quadratic forms of large discriminant, previous lower bound for the number of prime
curvatures is improved further in \hfill\break
[B-F2] to at least $\frac {cX}{(\log X)^{\frac 32 -\frac {\log 2} 2+\ve}}$.

The main result in this paper gives the correct order of magnitude.

\proclaim
{Theorem 2}
Given an integral primitive ACP, there is a positive $c$, such that for $X$ large the number of prime numbers less than $X$ which are curvatures of circles in the ACP, is
at least $c\frac X{\log X}$ (with $c>0$ an absolute constant). 
\endproclaim

Compared with the arguments due to Sarnak and refined in [B-F2], that are based on Iwaniec' theorem and representations by individual quadratic forms,
the strategy used here is different.
Our approach consists in introducing a generating function by considering the collected contribution of suitable families of binary quadratic forms
(constructed in \S0, \S1 of the paper).
These generating functions can then be analyzed using the circle method (in a rather standard way), to the extent of providing a main (arithmetical)
contribution with an error term.
In particular, we are able to establish Theorem 2 (relying also on the so-called `majorant property' for the set of the prime numbers, established in [B],
[G]).
The technique applied here may be organized better as to allow a treatment of the major-arcs contribution by spectral methods (using the
spectral analysis for the full Apollonian group), in the spirit of [B-K].
This leads to better error terms and statements that come close to a local to global principle.
That program is pursued in the forthcoming paper [B-K2].
Let us also mention the paper \cite {F-S} that gives evidence for the only congruence obstructions to appear $(\mod 24)$.

In this discussion, we should cite the paper [K-O], where counting results for the curvatures, with multiplicity, are obtained based on spectral
techniques (see also [BGS]).
In particular, it is shown in [K-O] that in any ACP $P$ the number of curvatures at most $X$ is of the order
$$
X^\delta \text { for } \ X\to\infty
$$
with $\delta =1,30068..$ is independent of the packing.
This amounts also to the number of quadruples bounded by $X$ in the orbit of the root quadruple under the Apollonian group $A$.

Let us briefly recall how binary quadratic forms enter the analysis (see [S]).
While $A$ is a `thin' (non-arithmetic) group, its subgroup $A_1=\langle S_2, S_3, S_4\rangle $ (= stabilizer of $x_1$) and similarly $A_2, A_3, A_4$ are
arithmetic.
More precisely, considering the map
$$
y=(y_2, y_3, y_4)= (x_2, x_3, x_4)+(a, a, a)
$$
the affine action of $A_1$ on $(x_2, x_3, x_4)$ is conjugated to the action of a finite index subgroup $\Gamma$ of $O_g(\Bbb Z)$, $g$ denoting the quadratic
form
$$
g(y)= y_2^2+y_3^2+y_4^2- 2y_2y_3 - 2y_2y_4 -2y_3y_4.
$$
By a further coordinate change
$$
A=y_2, B=\frac 12 (y_1-y_3+y_4), C=y_4
$$
transforming $g(y)$ in the quadratic form $\Delta(A, B, C)=B^2-AC$, $\Gamma$ is conjugated to the subgroup of $O_\Delta(\Bbb Z)$ generated by the reflections
$$
\left[ \matrix 1&-4&4\\ 0&-1&2\\ 0&0&1\endmatrix\right], \left[\matrix 1&0&0\\ 0&-1&0\\ 0&0&1\endmatrix\right], \left[\matrix 1&0&0\\ 2&-1&0\\ 4&-4&1
\endmatrix\right].
$$
Consider the spin double cover of $SO_\Delta(\Bbb Z)$ realized as image of $GL_2(\Bbb Z)$ under the homomorphism
$$
\rho:\pmatrix \alpha&\beta\\ \gamma&\delta\endpmatrix \mapsto \frac 1{\alpha\delta-\beta\gamma}
\left[\matrix \alpha^2& 2\alpha\gamma& \gamma^2\\
\alpha\beta& \alpha\delta+\beta\gamma &\gamma\delta\\
\beta^2& 2\beta\delta& \delta^2\endmatrix\right]
$$
with kernel $\pm I$. 
Then $\rho^{-1}(SO_\Delta (\Bbb Z)\cap \tilde\Gamma)$ contains $\left[\matrix 1&-2\\ 0&1\endmatrix\right], \left[\matrix
1&0\\ -2&1\endmatrix\right]$ and hence the principal congruence subgroup $\Lambda (2)$ of $SL_2(\Bbb Z)$.

It turns out that the set of values of $y_2=A$, $y_3=A+C-2B, y_4 =C$ contains at least those of $A$ with $(A, B, C)$ ranging in an orbit $\rho\big(SL_2(\Bbb
Z)\big) (A_0, B_0, C_0)^t$, i.e. the integers represented primitively by the binary quadratic form
$$
A_0 \alpha^2 +2B_0\alpha\gamma + C_0\gamma^2 \ \text { with } \ (\alpha, \gamma)= 1.
$$
The preceding provides an explicit recipe to produce curvatures in a given packing $P$.
Assume $(a_0, b_0, c_0, d_0)\in \Cal S=\Cal S(P)=A(a, b, c,d)^t$ and set
$$
A_0=a_0+b_0, 2B_0=a_0+b_0-c_0+d_0, C_0=a_0+d_0.
$$
Then all integers represented by the quadratic form
$$
A_0x^2+2B_0 xy+C_0y^2-a_0 \text { with } x, y \in\Bbb Z, (x, y)=1
$$
appears as curvatures of circles in the packing $P$.

This observation made in [S] plays a key role in [B-F1] and also in the construction of an appropriate 
family of binary quadratic forms described in \S0, \S1 of this paper.

\noindent
{\bf (0). Preliminary construction of a set of curvature quadruples}

Let $R_1$ be a large integer and denote $\Cal S_{R_1}$ the set of quadruples $(a, b, c, d) \in\Cal S=\Cal S(P)$ of the
Apollonian packing $P$ satisfying
$$
\max (|a|, |b|, |c|, |d|)\sim R_1.\tag 0.1
$$
Thus
$$
|\Cal S_{R_1}|> R_1^\delta \text { with } \ \delta > \frac {13}{10}.\tag 0.2
$$
Let
$$
R_2= R_1^{\frac 1{30}}.\tag 0.3
$$
Given $(a, b, c, d)\in\Cal S_{R_1}$, let $A= a+b, C=a+d, 2B=a+b-c+d$ and consider the set of integers
$$
S_{a, b, c, d}=\{ Ax^2 +2Bxy+Cy^2-a; x, y\in\Bbb Z, 0\leq x, y<R_2 \ \text { and }\ (x, y)=1\}.\tag 0.4
$$
Recall that $(A, B, C)= 1$ and $a^2=AC -B^2$ (by Descartes' equation).

As explained in [B-F1] the set $S_{a, b, c, d} $ is contained in the set of curvatures produced in the orbit of $(a, b, c, d)$ under group
elements of $A_1 =\langle S_2, S_3, S_4\rangle$ of  norm bounded by $R_2^2$.
Denote $\Cal S(a, b, c, d)\subset \Cal S_{R_1 R^2_2} \cap \langle S_2, S_3, S_4\rangle \pmatrix a\\ b\\ c\\ d\endpmatrix$ a set of quadruples $(a, b' , c', d')$ in
one-to-one correspondence with $S_{a, b, c, d}$ by projection on the $b'$-coordinate.

Thus for each $0\leq x, y<R_2$, $(x, y)=1$ there is some $g_{x, y}\in \langle S_2, S_3, S_4\rangle$ such that
$$
\Cal S(a, b, c, d)\subset \Big\{ g_{xy} \pmatrix a\\ b\\ c\\ d\endpmatrix; \ 0\leq x, y< R_2 \ \text { and } \ (x,y)=1\Big\}.\tag 0.5
$$
Obviously, if we fix $x, y$, all quadruples $g_{xy}\pmatrix a\\ b\\ c\\ d\endpmatrix$ are distinct and hence
$$
\sum_{(a, b, c, d)\in\Cal S_{R_1}} 1_{\Cal S(a, b, c, d)} \leq R_2^2.\tag 0.6
$$
Given $b'$, it follows from Descartes' equation that
$$
\pi_2^{-1} (b') \cap\Cal S_{R_1R_2^2}\ll R_1^{1+\ve} R^2_2\tag 0.7
$$
and (0.6), (0.7) imply that
$$
\sum_{(a, b, c, d) \in\Cal S_{R_1}} \ 1_{S_{a, b, c, d}} \ll R_1^{1+\ve} R^4_{2}.\tag 0.8
$$
To each $(a, b, c, d)\in\Cal S_{R_1}$, associate the distribution $\lambda_{a, b, c, d}$ on $\Bbb Z$ obtained as image
measure of 
$$
[0\leq x, y< R_2; (x, y) =1, f_a(x, y)\sim R_1 R_2^2 \ \text { and } \ \big(f_a(x, y), \prod_{p< R_2^{\frac 1{10}}} p)=1]
$$
under the map
$$
(x, y) \mapsto f_a(x, y) =Ax^2 +2Bxy+Cy^2 -a.
$$
Hence, supp\,$\lambda_{a, b, c, d}\in S_{a, b, c, d}$, $\Vert \lambda_{a, b, c, d}\Vert_\infty \ll R_1^\ve$ and elementary sieving
implies that certainly
$$
\frac {R_2^2}{(\log R_2)^2}\lesssim \Vert\lambda_{a, b, c, d}\Vert_1 \lesssim  R_2^2.\tag 0.9
$$
Also, from sieving, we obtain that for all $q\in \Bbb Z_+$
$$
\sum_{z\equiv  u (\mod q)} \lambda_{a, b, c, d}(z) \lesssim \Bigl(\frac 1q+\frac 1{R_2}\Bigr)^{\frac 1{10}} \Vert\lambda_{a, b, c,
d}\Vert_1\tag 0.10
$$
\big(in the argument, we distinguish the cases $q>(\log R_2)^{100}$ and $q< (\log R_2)^{100}$; we only need a crude estimate for our purpose\big).

Define
$$
\lambda=\sum_{(a, b, c, d)\in\Cal S_{R_1}} \ \lambda_{a, b, c, d}\tag 0.11
$$
which is a distribution on $[b' \in\Bbb Z; \, b'\sim R_1 R_2^2]$.
From (0.8), (0.9)
$$
|\Cal S_{R_1}| \frac {R_2^2}{(\log R_2)^2} <\Vert\lambda\Vert_1 \leq |\Cal S_{R_1}|. R_2^2  \ \text { and } \ \Vert\lambda\Vert_\infty
\ll R_1^{1+\ve} R_2^4.\tag 0.12
$$

By construction, $z\in \text{\,supp\,} \lambda$ has no prime factors less than $R_2^{\frac 1{10}}$ and from (0.10)
obviously
$$
\sum_{z\equiv u(\mod q)} \ \lambda(z) \lesssim \Bigl(\frac 1q+\frac 1{R_2}\Bigr)^{\frac 1{10}} \Vert\lambda\Vert_1.\tag 0.13
$$
Let $\eta_{a, b, c, d}$ be a distribution on $\Cal S(a, b, c, d)$ which image measure under  projection on the $b'$-coordinate
equals
$\lambda_{a, b, c, d}$ and set
$$
\eta=\sum_{(a, b, c, d)\in \Cal S_{R_1}} \ \eta_{a, b, c, d}.\tag 0.14
$$
Hence
$$
\Vert\eta\Vert_1 =\Vert\lambda\Vert_1.\tag 0.15
$$

Since clearly $\Vert\eta_{a, b, c, d}\Vert_\infty \leq\Vert\lambda_{a, b, c, d}\Vert_\infty \ll R_1^\ve$, it follows from (0.6)
that
$$
\Vert\eta\Vert_\infty \ll R_2^2 R_1^\ve\tag 0.16
$$
(0.13) may be rephrased as
$$
\sum_{b'\equiv u(\mod q)} \eta(a', b', c', d')\lesssim \Bigl(\frac 1q +\frac 1{R_2}\Bigr)^{\frac 1{10}}\Vert\eta\Vert_1
\text { for all } \ q\in\Bbb Z.\tag 0.17
$$

Next, we replace the distribution $\eta$ by a subset $\Cal C\subset\Cal S_{R_1R_2^2}$ which we
construct probabilistically by selecting $(a', b', c', d')\in\Cal C$ with probability
$$
\delta\eta(a', b', c', d')< 1\text { where } \ \delta=R_2^{-3} \quad \text {(cf. (0.16))}.
$$
By (0.12), (0.15), we obtain
$$
|\Cal C| \approx \delta. \Vert \eta\Vert_1 > R_2^{-2} |\Cal S_{R_1}|\tag 0.18
$$
and also
$$
\sum_{b'} |\Cal C_{b'}|^2 \lesssim \delta\Vert\eta\Vert_1 +\delta^2 \Vert\eta\Vert_1 \, \Vert\lambda\Vert_\infty
< |\Cal C|. R_1 R_2^2\tag 0.19
$$
with $\Cal C_{b'}$ denoting the fibers of $\Cal C$.

From (0.17) and standard large deviation inequalities, we deduce that for $q\in\Bbb Z_+$, $u\in\Bbb Z$
$$
\spreadlines{6pt}
\align
&|\{(a', b', c', d')\in\Cal C; b'\equiv u (\mod q)\}|\leq\\
&\delta \sum_{b'\equiv u(\mod q)} \ \eta(a', b', c', d') +c\sqrt{\log R_1} (\delta\Vert\eta\Vert_1)^{\frac 12}\lesssim\\
&\Bigl(\frac 1q+ \frac 1{R_2}\Bigr)^{\frac 1{10}}|\Cal C|.\tag 0.20
\endalign
$$
Relabeling the quadruples, we obtain a subset $\Cal C\subset\Cal S_R, R=R_1 R_2^2$ with the following properties

\roster
\item "{(0.21)}"  $|\Cal C| > R^{\delta-\frac 1{10}} > R^{6/5}$
\medskip

\item  "{(0.22)}" $a\sim R$ and $|b|, |c|,  |d|\lesssim R$ for $(a, b, c, d)\in\Cal C$
\medskip

\item  "{(0.23)}" $\sum_a|\Cal C_a|^2< R^{-1/5} |\Cal C|^2$
\medskip

\item  "{(0.24)}" For $(a, b,c,d)\in\Cal C$, $a$ has no prime factors less than $R^{\frac 1{320}}$ 
\medskip

\item  "{(0.25)}"
 $|\{(a, b, c, d)\in\Cal C; a\equiv u(\mod q)\}|\lesssim (q^{-1} +R^{-1})^{\frac 1{320}} |\Cal C| \ \text { for all }
\ q\in\Bbb Z_+$
\endroster

\bigskip

\noindent
{\bf (1). Introducing a family of quadratic forms}

Let $\Cal C\subset \Cal S_R$ be the set constructed in \S0.
Let $\Cal A=\pi_a (\Cal C)\subset \Bbb Z_+$.

To each $(a, b, c, d) \in\Cal C$, we associate again the binary form
$$
f(x, y)=Ax^2 +2Bxy+Cy^2\tag 1.1
$$
with
$$
A=a+b, C= a+d, 2B=a+b-c+d, \text { disc\,} f =- 4a^2, a^2 =AC-B^2.\tag 1.2
$$
Thus $(A, B, C)=1$.
Since $|A|, |B|, |C| \lesssim R$ and $|A.C|\gtrsim a^2\sim R^2$, it follows that $|A|, |C|\sim R$.

Denote $\Cal F$ the family of quadratic forms (1.1) obtained from $\Cal C$ and by $\Cal F_a\subset\Cal F$ those obtained from $\Cal C_a$.
We show that if we fix the discriminant, the number of equivalent forms in $\Cal F$ is $O(1)$.

Thus if $Ax^2 +2Bxy+Cy^2 $ and $A_1x^2+2B_1 xy+C_1y^2$ are equivalent, then
$$
\left\{
\aligned
&A_1= \alpha^2 A+2\alpha\gamma B+\gamma^2 C\\
&B_1= \alpha\beta A+(\alpha\delta+\beta\gamma)B+\gamma\delta C\\
&C_1 =\beta^2 A+2\beta\delta B+\delta^2 C
\endaligned
\right.
$$
for some $\pmatrix \alpha&\beta\\ \gamma&\delta\endpmatrix \in SL_2 (\Bbb Z)$.

Thus
$$
A_1 =A\Bigl(\alpha+\frac BA\gamma\Bigr)^2 +\frac {a^2}A\gamma^2.
$$
Hence
$$
\gamma^2 <\frac {AA_1}{a^2} <O(1)
$$
and since $|A|, |A_1|\sim R, \Big| \alpha+\frac BA\gamma\Bigr| <O(1), |\alpha|<O(1)+O(1)\frac {|B|}{|A|} < O (1)$.

Similarly, we see that, $|\beta|, |\delta|< O(1)$.
This shows that at most $O(1)$ quadratic forms obtained from $\Cal C_a$ are equivalent.
This proves our claim.

As a consequence, we obtain that for fixed $a$ and $M$
$$
\#\{(f, f_1, x, y, x_1, y_1)\in\Cal F_a\times\Cal F_a\times [1, M]^4; f(x, y) =f_1 (x_1, y_1)\}\ll (RM) ^\ve M^2|\Cal F_a|.\tag
1.3
$$
Indeed, fix $f\in \Cal F_a$ and $x, y$.
The integer $z= f(x, y)$ is at most $RM^2$ and is represented by at most $2^{\omega(z)}$ classes with discriminant $-4a^2$.
From the preceding, there are at most $O(1)$ \, $2^{\omega(z)}$ possibilities for $f_1\in\Cal F_a$ and for each $f_1, f_1(x_1, y_1)=z$
holds for at most $(RM)^\ve$ values of $(x_1, y_1)$.
This establishes (1.3).

\bigskip

\noindent
{\bf 2. Application of the circle method}

Let $R$ and $\Cal F$ be as in \S1.
Let $P$ be large and assume $ R<(\log P)^C$. 

For $f\in\Cal F$, let $\omega_f$ be the image measure on $\Bbb Z$ of $[1, P]^2$ under the map
$$
(x, y)\mapsto f_a (x, y)\ \text { with } \ f_a=f-a, -4 a^2=\text { disc\,} f
$$
under restriction $(x, y)=1$.

Actually it is technically more convenient to consider the image measure of \break
$1_{(x, y)=1} \gamma\big(\frac xP\big)\otimes\gamma \big(\frac yP\big)$,
where $0\leq\gamma\leq 1$ is a smooth bumpfunction supported  on $[0, 1]$.

Note that $\text{supp\,} \omega_f\subset [-RP^2, RP^2]$ (by construction of $\Cal F$) and we may assume \hfill\break
$\text{supp\,} \omega_f\subset
[0, RP^2]$ (by assuming $A>0$).
Define
$$
\omega =\underset {f\in \Cal F}\to {\overline\sum}\omega_f \ \text { where } \ \operatornamewithlimits{\overline\sum}_{f\in \Cal F} =\frac
1{|\Cal F|}
\operatornamewithlimits{\sum}_{f\in\Cal F}.
$$

The co-primality condition $(x, y)=1$ leads to technical complications (with no effect on the basic scheme of the argument).

Fix some integer $B=(\log P)^{10}$ and replace the restriction $1_{(x, y)=1}$ by
$$
\sum_{d|(x, y), d< B} \mu(d).
$$
This expression equals 1 if $(x, y) =1$, vanishes if $1<(x, y)<B$ and is bounded by $\tau\big((x, y)\big)$.
Hence this replacement introduces an error at most 
$\sum_{x, y<P} \tau \big((x.y)\big) 1_{(x, y)>B} \ll \frac{P^2}B(\log P)^2$ for the counting functions $\omega_f$ and
$\omega $ (in $\ell^1(\Bbb Z)$-norm), which is harmless.

With the above modification, we obtain
$$
\hat\omega(\theta) =\sum_{d<B} \mu (d)\Big\{\operatornamewithlimits{\overline\sum}
\limits_{\Cal F}\sum_{d|(x, y)} \gamma \Big(\frac xP\Big)\gamma\Big(\frac yP\Big) e\big(f_a(x, y)\theta)\Big\}
$$
and 
$$
|\hat\omega(\theta)|^2 \ll \sum_{d< B} d^{1+\ve}\Big|\operatornamewithlimits{\overline\sum}\limits
_{\Cal F} \sum_{d|(x, y)} \gamma\Big(\frac xp\Big)\gamma\Big(\frac yP\Big) e\big(f_a(x, y)\theta\big)
\Big|^2.\tag 2.1
$$

Fixing $d$ square free, we have to analyze the expressions
$$
S_\omega(\theta) =\operatornamewithlimits{\overline\sum}_{\Cal F} \sum_{d|(x, y)}  \, \gamma\Big(\frac xP\Big)\gamma\Big(\frac yP\Big)
e\big((Ax^2+2Bxy+Cy^2-a)\theta\big).\tag 2.2
$$
For $q<RP$, denote for $(q, b)=1$
$$
\Cal M(q, b)=\Big[\Big|\theta-\frac bq\Big| < \frac 1{qRP}\Big] \subset\Bbb T=\Bbb R/\Bbb Z.
$$
Our main concern is to obtain suitable bounds on
$$
\sum_{q\sim Q} \sum_{(b, q)=1}  \ \int\limits_{\Cal M(q, b)} |S_\omega (\theta)|^2 d\theta.\tag 2.3
$$
Let $\theta =\frac bq+\vp\in \Cal M(q, b), |\vp|<\frac 1{qRP}$.
Then
$$
f_a (x, y) \theta= f_a(x, y) \frac bq+ f_a(x, y)\vp
$$ 
and
$$
\spreadlines{6pt}
\align
&S_{\omega_f}(\theta)\equiv \sum_{d|(x, y)} \gamma \Big(\frac xP\Big)\gamma \Big(\frac yP\Big) \ e\big(f_a(x, y)\theta\big)=\\
&\sum_{\Sb 0\leq k, \ell< q\\ d_0|(k, \ell)\endSb} e\Big(f_a(k, \ell)\frac bq\Big)\Big[\sum_{{\matrix x\equiv k, y\equiv \ell (\mod q)\\
 x\equiv y\equiv 0 (\mod d_1)\endmatrix}} e\big(f_a(x, y) \vp\big) \gamma\Big(\frac xP\Big) \gamma\Big(\frac yP\Big)\Big]\tag 2.4
\endalign
$$
where $d_0=(d, q), d = d_0d_1$.

We distinguish 2 cases.
\medskip

\noindent
{\bf Case I.}  \ $q<P$.

Rewrite the second factor in (2.4) as
$$
\sum_{\Sb r, s\lesssim\frac PQ\\ rq+k\equiv sq+\ell\equiv 0(\mod d_1)\endSb} \gamma\Big(\frac {rq+k}P\Big)\gamma\Big(\frac {sq+\ell}P\Big)
e\big(f(rq+k, sq+\ell)\vp\big)
$$
(dropping a multiplicative factor)
$$
=\sum_{r', s' \lesssim \frac P{d_1 Q}} \gamma\Big(\frac{r'd_1q+k'}P\Big)\gamma\Big(\frac {s' d_1q+\ell'} P \Big) e\big(f(r'd_1 q+k', s' d_1 q+\ell')\vp\big)
$$
with $0\leq k', \ell'< d_1q; k'\equiv k(\mod q), \ell'\equiv \ell(\mod q), k'\equiv\ell' \equiv 0 (\mod d_1)$.

From the Poisson summation formula, we obtain
$$
\align
&\frac 1{q^2d_1^2} \sum_{m, n\in\Bbb Z}\iint \gamma\Big(\frac {y+k'}P\Big) \gamma\Big(\frac {z+\ell'}P\big) e\big(f(y+k', z+\ell')\vp\big)
e\Big(m\frac y{qd_1}+n\frac z{qd_1}\Big) dydz\\
&=\frac 1{q^2d_1^2} \sum_{m, n\in\Bbb Z} J_f(qd_1, m, n, \vp) e_q(-m\td_1 k - n\td_1\ell)\tag 2.5
\endalign
$$
where $d_1\td_1 \equiv 1 (\mod q)$ and
$$
J_f(q, m, n,\vp)=\iint \gamma\big(\frac yP\Big)\gamma\Big(\frac zP\Big) e\big(f(y, z)\vp\big) e\Big(\frac mq y+\frac nqz\Big)dydz.\tag 2.6
$$
Note that, by stationary phase
$$
|J_f(q, m, n, \vp)|\lesssim \frac 1{|\vp| \, |\text{discr } (f)|^{1/2}}\lesssim \frac 1{R|\vp|}.\tag 2.7
$$
Also, in (2.5) the significant contributions come from values $m ,n$ satisfying
$$
\frac {|m|}{qd_1}, \frac {|n|}{qd_1}\lesssim |\nabla f|. |\vp|\lesssim R.P. \frac 1{QRP}\tag 2.8
$$
hence $|m|, |n|=0(d_1)$.

From (2.7), there is an obvious bound $0\big(\frac 1{q^2} \frac 1{R|\vp|}\big)$ on (2.5) and $0\big(\frac 1{R|\vp|}\big)$ on (2.4), (2.2).
Substituting in (2.3), we see that the contribution of $|\vp|>\vp_*$ is at most
$$
\frac{B^2Q^2}{R^2\vp_*} < \frac {P^2}R (\log P)^{-C} \ \text { provided } \ \vp_*> (\log P)^C\frac {Q^2}{P^2}.
$$
If $Q<P^{\frac 1{10}}$, we may therefore restrict $|\vp|< P^{-\frac 53}$ and obtain $m=n=0$ in (2.8).  Thus we distinguish the ranges

\noindent
$\underline{Q<P^{\frac 1{10}}}$.
Then (2.4) is replaced by
$$
\frac 1{q^2d^2_1} \Big[\sum_{\Sb 0\leq k, \ell<q\\ d_0|(k, \ell)\endSb} e\Big(f_a(k, \ell)\frac bq\Big)\Big] J_f(\vp)\tag 2.9
$$
contributing in (2.3) for
$$
\frac 1{Q^4} \sum_{q\sim Q} \frac 1{d_1^4} \sum_{(b, q)=1} \Big|\operatornamewithlimits{\overline\sum}\limits
_{\Cal F} c_f\sum_{\Sb 0\leq k, \ell<q\\ d_0|(k, \ell)\endSb}
e\Big(f_a(k, \ell)\frac bq\Big)\Big|^2. \int_{|\vp|<\frac 1{RPQ}}\Big(\min \Big(P^2, \frac 1{R|\vp|} \Big)\Big)^2 d\vp\tag 2.10
$$
(where $|c_f|\leq 1$).

Summing over $d$, we obtain
$$
\frac {P^2}R\sum_{d_0} d_0^{1+\ve} \Big\{\sum_{\Sb q\sim Q\\ d_0|q\endSb} \sum_{(b, q)=1} \Big|
\operatornamewithlimits{\overline\sum}\limits_{\Cal F} c_f S_{f, d_0} (q, b, 0, 0)\Big|^2\Big\}
\tag 2.11
$$
denoting
$$
\align
S_f(q, b, m , n)&=\frac 1{q^2} \sum_{0\leq k, \ell<q} e_q\big(bf_a(k, \ell)-mk-n\ell\big)\tag 2.12\\
S_{f, d_0} (q, b, m, n)&=\frac 1{q^2} \sum_{\Sb 0\leq k, \ell < q\\ d_0|(k, \ell)\endSb} e_q \big(b f_a(k, \ell)-mk-n\ell\big)\tag 2.12'
\endalign
$$
and where $d_0<B$ is a square-free divisor of $q$.
\medskip

\noindent
$\underline {P^{\frac 1{10}}\leq Q\leq P}$

By (2.5), (2.4) becomes
$$
\frac 1{d_1^2} \sum_{|m|, |n|<0(d_1)} S_{f, d_0}(q, b, m\td_1, n\td_1) \, J_f(qd_1, m, n, \vp)\tag 2.13
$$
and the contribution in (2.3) is at most
$$
B^{2+\ve} \, \frac {P^2}{R} \max_{d_0, d_1, m, n} \sum_{\Sb q\sim Q\\ d_0|q\endSb}\Big\{\sum_{(b, q)=1} \Big|
\operatornamewithlimits{\overline\sum}\limits_{\Cal F} c_{f, q} S_{f, d_0} (q, b, m\td_1,
n\td_1)\Big|^2\Big\}\tag 2.14
$$
with $|c_{f, q}|\leq 1, d, d_1$ square free, $(d_0, d_1)=1, d_0d_1 <B$ and $m, n=0(d_1)$.

\bigskip
{\bf Case II. \  $P<Q<PR$}
\medskip

Since for $\theta =\frac bq+ \vp \in\Cal M(q, b)$
$$
f_a (x, y) \theta= f_a(x, y) \frac bq + f_a(x, y)\vp \ \text { with } \ |f_a(x, y) | \, |\vp| <R.P^2 .
\frac 1{RP Q} < O(1)
$$
we may replace $\theta$ by $\frac bq $ (dropping $\vp$) and $S_\omega (\theta)$ becomes
$$
\operatornamewithlimits{\overline\sum}_{\Cal F}\sum_{\Sb 1\leq k, \ell \leq P\\ d|(x, \ell)\endSb} 
e_q \big(bf_a(k, \ell)\big).\gamma\Big(\frac kP\Big)\gamma\Big(\frac\ell P\Big).\tag 2.15
$$

Thus the inner sum in (2.15) equals
$$
\sum_{d_0|(k, \ell)} e_q \big(bf_a(d_1 k, d_1\ell)\big) \gamma\Big(\frac {d_1 k}P\Big)\gamma\Big(\frac {d_1\ell}P\Big)
$$
and completing the sum, we obtain
$$
\align
&e_q(-ab).\operatornamewithlimits{\overline\sum}\limits
_{|u|, |v|<\frac QP d_1} \sum_{\Sb 0\leq x, y< q\\d_0|(x, y)\endSb} e_q\big(d_1^2 bf(x, y)+ux+vy)\\
& = q^2\operatornamewithlimits{\overline\sum}\limits_{|u|, |v|<\frac QP d_1} S_{f, d_0} (q, b, \td_1 u, \td_1 v).\tag 2.16
\endalign
$$
The contribution to (2.3) may be bounded by
$$
\frac {B^{2+\ve} Q^3}{RP} \ \sum_{\Sb q\sim Q\\ d_0|q\endSb} \ \sum_{(b, q)=1} \Big|
\operatornamewithlimits{\overline\sum}\limits_{u, v<\frac QP d_1} \operatornamewithlimits{\overline \sum}\limits_{\Cal F} S_{f, d_0}
(q, b, \td_1 u, \td_1 v)\Big|^2\tag 2.17
$$
with $d_0, d_1$ square free, $(d_0, d_1)=1$ and $d_0d_1<B$.

\bigskip

\noindent
{\bf (3). Evaluation of the Gauss sum}

Analyzing further (2.11), (2.14), (2.17), we evaluate expressions of the form
$$
\sum_{(b, q)=1} S_f (q,b, u,v)  \overline{S_{f_1}(q, b, u_1,v_1)} \tag 3.1
$$
and also
$$
\sum_{(b, q) =1} S_{f, d_0} (q, b, u, v) \overline{S_{f_1, d_0}(q, b, u_1, v_1)}.
\tag 3.1$'$
$$
Consider first (3.1)

Factoring $q$ as a product of prime power  $p^r$, (3.1) factors correspondingly.
Recall that
$$
S_f(q, b, u, v)=\frac 1{q^2}\sum_{0\leq x, y<q} e_q (bf(x, y)+ux+vy-ba)\tag 3.2
$$
with $f(x, y) =Ax^2+2Bxy+Cy^2, (A, B, C)= 1$.

Let $q=p^r$ be a prime power.  We may assume $(A, p)=1$.

Write (assuming $p\not =2$; for $p=2$ there are some extra technicalities that we omit here and to which we will return in \S7)
and using the notation $\overline { \cdot  }$ for the multiplicative inverse $(\mod q)$
$$
\spreadlines{6pt}
\align
&bf(x, y) +ux+vy-ab= b(Ax^2+2Bxy+Cy^2)+ux+vy-ab\\
&= bA(x+ B\bar A\, y)^2+b a^2\bar A \, y^2 +ux+vy-ab\\
&=bA(x+B\bar A \, y+ \overline{2bA} \, u)^2 +ba^2 \bar A \, y^2+(v-B\bar A \, u) y- \overline{4bA} \, u^2-ab \quad (\mod q).
\endalign
$$
Hence, by Gauss sum evaluation, we obtain (cf. [BEW])
$$
S_f(q, b, u, v) \sim\frac 1{q^{3/2}} \Big(\frac{bA}{q}\Big) e_q (-\overline{4Ab} \, u^2-ab)\Big[\sum_{0\leq y< q} e_q (ba^2\bar A
\, y^2+
(v-B\bar A \, u) y\big)\Big].\tag 3.3
$$
Let $\frac {a^2} q=\frac {\tilde a}{\tilde q}$ with $(\tilde a, \tilde q)= 1$.
Thus writing $y =\tilde y+z\tilde q, 0\leq \tilde y<\tilde q$, $0\leq z\leq\frac q{\tilde q} = (q, a^2)$
$$
\sum_{0\leq y<q} \cdots = \sum_{0< \tilde y<\tilde q} \ \sum_{0\leq z< (a^2, q)} e\Big(b\bar A \, \frac{\tilde a(\tilde y)^2}{\tilde q}
+
(v-B\bar A \, u)\frac {\tilde y} q +(v-B\bar A \, u)\frac z{(a^2, q)}\Big)\tag 3.4
$$
and $(3.4)=0$ unless  $(a^2, q)| v-B\bar A u$, in which case we set 
$$
v-B\bar A u= (v-B\bar A u)^\sim . (a^2, q)
$$
Thus (3.4) becomes
$$
\spreadlines{6pt}
\align
&(a^2, q)  .\sum_{0\leq y<\tilde q} e_{\tilde q} \big(b\bar A \tilde a y^2 +(v-B\bar A u)^\sim y)\\
& = (a^2, q) \sum_{0\leq y<\tilde q} e_{\tilde q} \big ( b\bar A\tilde a(y+ (Av-Bu)^\sim \, \overline{2b\tilde a} )^2 -\overline{4
b\tilde a} \,
\bar A \big((Av-Bu)^\sim \big)^2\big)\\
&\sim (a^2, q) (\tilde q)^{1/2} \Big(\frac {b\bar A\tilde a}{\tilde q}\Big) e_{\tilde q} \big(-\overline {4b\tilde aA} \,
\big((Av-Bu)^\sim)^2\big).\tag 3.5
\endalign
$$
Hence, from (3.3), (3.5)
$$
\spreadlines{6pt}
\align
S_f (q, b, u, v)&=\frac {(a^2, q)^{1/2}}q \ \Big(\frac {bA} q\big)\Big(\frac {b\bar A\tilde a}{\tilde q}\Big) e_q
(-\overline {4Ab} \, u^2 -ab) e_{\tilde q} (-\overline{4b\tilde a A}\, \big(( Av-Bu)^\sim\big)^2\big)\\
&=\left\{\aligned \frac {(a^2, q)^{\frac 12}}q \ e_q(-\overline{4Ab} \, u^2-ab) e_{\tilde q} (-\overline{4 b\tilde a A} \,
\big((Av-Bu)^\sim \big)^2\big) \ & \ \text { if }
\tilde q\not =1\\
\frac{(a^2, q)^{\frac 12}}q \Big( \frac{bA} q\Big) e_q \big(-\overline {4Ab} \, u^2 -ab)\qquad\qquad\qquad\qquad & \ \text { if }
q|a^2.\endaligned\right.\tag 3.6
\endalign
$$

It follows that
$$
(3.1) =\frac {(a^2, q)^{\frac 12} (a^2_1, q)^{\frac 12}} {q^2}\sum_{(b, q)=1} e_q \Big((\overline{4A_1} \, u_1^2 -\overline{4A} \,
u^2)\bar b + (a_1 -a) b\Big)
\overline {E(b)} \, E_1(b)\tag 3.7
$$
where $E(b)$ is either
$$
\Big(\frac{bA}q\Big) \text { or } \ e_{\tilde q} \big(-\overline{4A\tilde a} \, \big((Av-Bu)^\sim\big)^2\bar b\big) \text
{ with } \tilde q=\frac q{(q, a^2)}
\tag 3.8
$$
where we assume $(q, a^2)|Av-Bu$ and similarly for  $E_1(b)$.
\medskip

Thus the sum in (3.7) is of Kloosterman or Sali\'e-type.
\medskip

There is the following elementary estimate (which suffices for our needs)
$$
\Big|\sum_{\Sb (b, p)=1\\ 0<b< p^r\endSb} e_{p^r} (cb+d \bar b)\Big|< (p^r)^{3/4} (p^r, c, d)^{1/4}\tag 3.9
$$
and also
$$
\Big|\sum_{\Sb (b, p)=1\\ 0<b<p^r\endSb}        \Big(\frac bp\Big) e_{p^r} (cb+d\bar b)\Big| < (p^r)^{3/4} (p^r, c, d)^{1/4}.\tag
3.10
$$
Hence we may state the following bound for $q=p^r$
$$
(3.7) < (a^2, q)^{\frac 12} (a_1^2, q)^{\frac 12} (q, a-a_1)^{1/4} q^{-5/4}\tag 3.11
$$
and also, for $a=a_1$ 
$$
\spreadlines{2pt}
\align
&(3.7) < \\
& q^{-5/4} (a^2, q) \bigl(q, \overline {4A_1} u^2_1 -\overline{4A} u^2
+(q, a^2)\overline {4A_1\tilde a} \big((A_1v_1-B_1u_1)^\sim\big)^2
-(q, a^2)\overline{4A\tilde a} \big((Av-Bu)^\sim\big)^2 \big)^{1/4}\tag 3.12
\endalign
$$
where $(a^2, q) |Av -Bu, (a^2, q) |A_1v_1 -B_1u_1$ and $\tilde a =\frac {a^2}{(a^2, q)}, (Av-Bu)^\sim =\frac {Av-Bu}{(a^2, q)}$.

Consider next (3.1$'$) for which there is again factorization.

If $q=p^r$ and $p|d_0$, then
$$
S_{f, p}(p^r, b, u, v) = p^{-2r} e_{p^r} (-ba) \ \sum_{0\leq x, y<p^{r-1}} e_{p^{r-1}}\big(bf(x, y)p+ux+vy).\tag 3.13
$$
Note that by assumption (0.24) and since $p<(\log P)^{10}$, $(a, p)=1$.

We distinguish several cases

\noindent
(3.14) $\underline{r=1}$

Then $|(3.13)|= p^{-2}$ and $|(3.1')|= (p-1) p^{-4}$, which is certainly bounded by
$$
p^{-\frac 74}  \ (3.11)\tag 3.15
$$
and
$$
p^{-\frac 74} \ (3.12).\tag 3.16
$$

\noindent
(3.17) $\underline {r\geq 2}$

Clearly (3.13) vanishes unless $p|(u, v)$.
Writing $u=pu_1, v= pv_1$,
$$
(3.13) = e_{p^r}(-ba) p^{-2} p^{-2(r-2)} \Big[\sum_{0\leq x, y< p^{r-2}} e_{p^{r-2}}(bf(x, y) +u_1x+v_1y)\Big].\tag 3.18
$$
If $r>2$, repeating the analysis of the exponential sum with $q$ replaced by $q_1=p^{-2} q$, we obtain instead of (3.6)
$$
\frac 1q e_{q_1} \big(-\overline {4Ab} u^2 -\overline{4ba^2A} \big((Av-Bu)^\sim \big)^2\big) e_q(-ab).\tag 3.19
$$
This gives for (3.1) the bound
$$
(q, a-a_1)^{1/4} q^{-5/4}.\tag 3.20
$$
Factor
$$
q=\prod_{(p, d_0)=1} p^r \cdot \prod_{p|d_0} p\cdot \prod_{\Sb p|d_0\\ r>1\endSb} p^r
$$
and consider the corresponding factorization of (3.1').
Apply (3.11), (3.12) if $(d, p)=1$, (3.15), (3.16) if $p|d_0, r=1$ and (3.20) if $p|d_0, r>1$.

\bigskip

\noindent
{\bf (4). Estimation of (2.11), (2.14)}

Expressing the square of the inner sum and carrying out the summation over $b$, we evaluate  (3.1), (3.1$'$) by (3.11), (3.15), (3.20).

Hence for $\sum_{(b, q)=1} |\cdots|^2 $ we obtain the bound
$$
Q^{-5/4} \operatornamewithlimits{\overline\sum}\limits_{f, f_1\in\Cal F} (a^2, q)^{1/2} (a_1^2, q)^{1/2} 
(q, a-a_1)^{1/4} \Big(\prod_{\Sb p|d_0\\ r=1\endSb} p^{-7/4}\Big).\tag 4.1
$$
Substitution in (2.11) gives
$$
\align
&\frac {P^2}{RQ^{5/4}} \operatornamewithlimits{\overline\sum}\limits_{f, f_1\in\Cal F} 
\sum_{d_0} \sum_{\Sb q\sim Q\\ d_0|q\endSb} (a^2, q)^{\frac 12} (a_1^2, q)^{\frac 12}
(a-a_1, q)^{\frac 14}
\Big(\prod_{\Sb p|d_0\\r=1\endSb} p^{-\frac 34+\ve}\Big) \Big(\prod_{\Sb p|d_0\\r>1\endSb} p^{1+\ve}\Big)\\
& \leq \frac {P^2}{RQ^{5/4}} \operatornamewithlimits{\overline\sum}\limits_{f, f_1 \in\Cal F} 
\sum_{\Sb q\sim Q\\ d_0|q\endSb} (a^2, q_1)^{\frac 12} (a^2_1, q_1)^{\frac 12}
(a-a_1, q_1)^{\frac 14} d_{0, 2}^{\frac 54+\ve}\tag 4.2
\endalign
$$
where $q=q_1 d_{0, 2}$ and $d_{0, 2} =\prod_{\Sb p|d_0\\r>1\endSb}p$.
Note that $d_0|q_1$ and hence for given $q_1$, there
are at most $\min (q_1^\ve, B^{1+\ve})$ possibilities for $d_{0, 1}, d_{0, 2}$ and $q$.

Specify $q_1\sim Q_1$ where $Q_1<Q<\min (BQ_1, Q_1^2)$ and replace (4.2) by
$$
\frac {P^2}{RQ_1^{5/4}} \min (B^2, Q_1^\ve) \operatornamewithlimits{\overline\sum}_{f, f_1 \in \Cal F} \sum_{q_1\sim Q_1} (a^2, q_1)^{\frac 12} 
(a^2_1, q_1)^{\frac 12} (a-a_1, q_1)^{\frac 14}.\tag 4.3
$$
\medskip

Next, we evaluate
$$
\frac {P^2}{RQ^{5/4}} \ \operatornamewithlimits{\overline\sum}_{f, f_1\in\Cal F} \ \sum_{q\sim Q} (a^2, q)^{1/2} (a^2_1, q)^{1/2} (a-a_1, q)^{1/4}.\tag 4.4
$$
The contribution for $a=a_1$ may be bounded by
$$
\align
&\frac {P^2}{RQ} \ \frac 1{|\Cal F|^2} \ \sum_{a\in\Cal A} |\Cal F_a|^2\sum_{q\sim Q} (a^2, q)\leq\\
&\frac {P^2}{RQ} \ \frac 1{|\Cal F|^2} \ \sum_{a\in\Cal A} |\Cal F_a|^2 \tau (a^2)Q\ll \frac {P^2}{R^{1-\ve}} \
\frac 1{|\Cal F|^2} \sum_a |\Cal F_a|^2 \overset{(0.23)}\to \ll \frac {P^2}R R^{-\frac 15+\ve}.\tag 4.5
\endalign
$$

For $a \not= a_1$, we obtain
$$
\frac {P^2}{RQ^{5/4}}\ \frac 1{|\Cal F|^2} \sum_{\Sb a, a_1\in\Cal A\\ a\not= a_1\endSb} \ |\Cal F_a| \, |\Cal F_{a_1} |
\sum_{q\sim Q} (a^2, q)^{1/2} (a_1^2, q)^{1/2} (a-a_1, q)^{1/4}.\tag 4.6
$$
Clearly, for fixed $a\not= a_1$
$$
\sum_{q\sim Q} (a^2, q)^{1/2} (a_1^2, q)^{1/2} (a-a_1, q)^{1/4} \leq \sum_{\Sb d|a^2, d_1|a_1^2\\ [d, d_1]\lesssim Q\endSb}
d^{1/2} d_1^{1/2} \frac {Q^{5/4}}{[d, d_1]} \lesssim Q^{5/4}\tag 4.7
$$
since $a, a_1$ are pseudo-prime by (0.24).

Also
$$
\spreadlines{6pt}
\align
\sum_{q\sim Q} (a^2, q)^{1/2} (a_1^2, q)^{1/2} (a-a_1, q)^{1/4}&\leq \sum_{q\sim Q} (a, q) (a_1, q) (a-a_1, q)^{1/4}\\
&\lesssim \sum_{\Sb d| a, d_1| a_1, d'|a-a_1\\ [d, d_1, d']\lesssim Q\endSb} d d_1(d')^{1/4} \frac Q{[d, d_1, d']}.\tag 4.8
\endalign
$$
Denote $f=(d, d_1)$.
Thus $f|d'$ and since $f^{-1} d, f^{-1} d_1, f^{-1} d'$ are pairwise coprime, $[d, d_1, d']\geq f^{-3} dd_1 d'$ and therefore
$$
(4.8)\lesssim Q.(a, a_1)^3 \sum_{\Sb d'|a-a_1\\ d'\lesssim Q\endSb} (d')^{-3/4}.\tag 4.9
$$

From (4.7), (4.9)
$$
(4.6)\lesssim \frac {P^2} R \ \frac 1{|\Cal F|^2} \sum_{\Sb a,  a_1 \in \Cal A\\ a\not= a_1\endSb} |\Cal F_a|.|\Cal F_{a_1}| \min 
\Big(1, Q^{-\frac 14} (a, a_1)^3 \sum_{\Sb d'|a-a_1\\ d'<Q\endSb} \Big(\frac 1{d'}\Big)^{\frac 34}\Big).\tag 4.10
$$
We distinguish two cases.

Assume $(a, a_1)=\Delta>Q^{10^{-4}}$.
Estimate
$$
(4.10)\lesssim \frac {P^2}R\ \frac 1{|\Cal F|^2} \sum_{a\in\Cal A} |\Cal F_a|.\sum_{\Delta|a, \Delta>Q^{10^{-4}}} \ \sum_{a_1\equiv a(\mod\Delta)}
|\Cal F_{a_1}|.\tag 4.11
$$

Again, since $a$ is pseudo-prime, $\Delta$ is restricted to $0(1)$ values, once $a$ fixed.
From (0.25),
$$
|\{ f\in\Cal F; a\equiv u(\mod \Delta)\}|\lesssim \Delta^{-\frac 1{320}}|\Cal F|
$$
and hence
$$
(4.11) < \frac {P^2}R Q^{-10^{-7}}.\tag 4.12
$$
Assume $(a, a_1)\leq Q^{10^{-4}}$. Then
$$
\align
(4.10)&\lesssim \frac {P^2}R Q^{-\frac 14 +3.10^{-4}} \ \frac 1{|\Cal F|^2} \sum_{\Sb a, a_1 \in\Cal A\\ a\not= a_1\endSb} |\Cal F_a| \ |\Cal F_{a_1}|
\sum_{\Sb d'|a-a_1\\ d'< Q\endSb} \Big(\frac 1{d'}\Big)^{3/4}\\
&\lesssim \frac {P^2}R Q^{-\frac 14 +3.10^{-4}} \sum_{d'<Q} \big(\frac 1{d'}\Big)^{3/4} \Big(\frac 1{d'}\Big)^{\frac 1{320}}\\
&< \frac {P^2}R Q^{-\frac 1{400}}.\tag 4.13
\endalign
$$
From (4.5), (4.12), (4.13), we obtain
$$
(4.4) <\frac {P^2}R \big(R^{-\frac 16}+Q^{-10^{-7}}\big).\tag 4.14
$$

We assume $R\sim (\log P)^C$ with $C$  a sufficiently large constant.

Since $B\sim (\log P)^{10}$, it follows from (4.14) that
$$
(4.3)<\frac {P^2}{R} \big(R^{-\frac 17}+Q^{-10^{-8}}\big).\tag 4.15
$$
The same bound also holds for (2.14).

\bigskip
\noindent
{\bf (5). Estimation of (2.17)}

Expressing the square of the inner sums in (2.17) gives
$$
\frac {B^{2+\ve} Q^3}{RP} \sum_{\Sb q\sim Q\\ d_0|q\endSb} \ \operatornamewithlimits{\overline\sum}\limits_{f, f_1\in\Cal F} \ 
\operatornamewithlimits{\overline\sum}\limits _{u, v, u_1, v_1< \frac QP d_1}
\ \sum_{(b, q)=1} S_{f, d_0} (q, b, \td_1 u, \td_1 v)\overline {S_{f_1, d_0} (q, b, \td_1 u_1, \td_1 v_1)}\tag 5.1
$$
where the inner sum is of type (3.1$'$).

Note that here $d_0, d_1$ are fixed.

To estimate the contribution for $a\not=a_1$, use again the bound (3.11) on (3.1) (which from previous discussion is always valid).
We obtain
$$
\spreadlines{6pt}
\align
B^{2+\ve} . \frac {Q^{7/4}}{RP} &\ \operatornamewithlimits{\overline\sum}_{f, f_1\in\Cal F} \ \sum_{q\sim Q} (a^2, q)^{\frac 12}
(a^2_1, q) ^{\frac 12} (a-a_1, q)^{1/4}\\
&\ll \frac 1P Q^{\frac {11}4} R^{\frac 54} B^{2+\ve} < P^{7/4} R^5 <\frac {P^2} R R^{-1}.\tag 5.2
\endalign
$$
Next the $a=a_1$ contribution

Let $q=q_1q_2$ with $q_2=\prod_{\Sb p|d_0\\r>1\endSb} p^r$ and factor (3.1$'$) according to $q=\prod_p^r$.
If $p|q_1$, the bound (3.12) applies to (3.1$'$) with $q=p^r$.

For $p|q_2$, apply (3.11) which gives the bound $q_2^{-1}$ on the $q_2$-factor.
Let $q_1\sim Q_1, q_2\sim Q_2$ (noting that, by construction, the number of $q_2$-values is at most $Q_2^{\frac 12}$).
Thus the contribution to (5.1) may be bounded by
$$
\frac {B^{2+\ve} Q_1^3 Q_2^{5/2}}{RP} \sum_{\Sb q_1\sim Q_1\\ d_{0, 1}|q_1\endSb} \frac 1{|\Cal F|^2} \sum_{a\in\Cal A} \ \sum_{f, f_1\in\Cal F_a} \
\operatornamewithlimits{\overline\sum}\limits_{u, v, u_1, v_1<\frac QP d_1} \  (5.3)\tag 5.4
$$
with
$$
(5.3) =\prod_{p|q_1} \sum_{(b, p)=1} S_{f, d_{0, 1}}(p^r, b, \td_1 u, \td_1 v)\overline{S_{f_1, d_{0, 1}} (p^r, b, \td_1 u_1, \td_1v_1)}
$$
and $d_{0, 1} =\prod_{p|d_0, r=1} p$.

Because $(d_0, d_1)=1$, the factor $\td_1$ in (5.3) turns out to be irrelevant and we drop it for simplicity.

Applying (3.12) for each prime $p|q$, we obtain
$$
\align
&|(5.3)|< (a^2, q_1) Q_1^{-5/4} \prod_{p| q_1} \\
&\big( p^r, \overline{4A_1}\, u_1^2  -\overline{4A} u^2+(p^r, a^2) \overline{4A_1 \tilde a}
\big((A_1 v_1-B_1 u_1)^\sim \big)^2- (p^r, a^2) \overline{4A\tilde a} \big((Av-Bu)^\sim)^2\big)^{1/4}.\tag 5.5
\endalign
$$
We distinguish several cases for the factors in (5.5).

Assume
$$
 \big(p^r, \overline{4A_1}\, u_1^2 -\overline{4A} u^2 +(p^r, a^2) \overline{4A_1\tilde a} \big((A_1v_1-B_1u_1)^\sim\big)^2
-(p^r, a^2) \overline{4A\tilde a} \big((Av-Bu)^\sim\big)^2\big)>p^{[\frac r2]}.\tag 5.6
$$
Multiplying with $a^2$ gives then
$$
\overline{4A} u^2a^2+\tilde a \, \overline{4A\tilde a} (Av -Bu)^2 \equiv \overline{4A_1} u_1^2 a^2+\tilde a \, \overline{4A_1\tilde a}
(A_1v_1 -B_1u_1)^2 \quad (\mod p^{[\frac r2]+1})
$$
and since $\tilde a\bar{\tilde a}\equiv 1(\mod\frac {p^r}{(a^2, p^r)})$ and $Av-Bu \equiv A_1v_1 -B_1u_1\equiv 0 \quad\big(\mod (a^2,
p^r)\big)$
$$
\bar A u^2 a^2 +\bar A(Av -Bu)^2 \equiv \bar A_1 u_1^2 a^2 +\bar A_1(A_1v_1 -B_1 u_1)^2 \quad (\mod p^{[\frac r2]+1}).\tag 5.7
$$
Hence, since $a^2=AC-B^2 =A_1C_1 -B_1^2$,
$$
f(v, -u)\equiv f_1(v_1, -u_1) \ (\mod p^{[\frac r2]+1}). \tag 5.8
$$
Since $|A|, |B|, |C|, |u|, |v|\lesssim R$, it follows that either
$$
f(v, -u)=f_1(v_1, -u)\tag 5.9
$$
or (5.8) can only hold for a set $\sigma$ of primes $p|q$ such that
$$
\prod_{p\in\sigma} p^{[\frac r2]+1} < R^4.\tag 5.10
$$
Hence, if (5.9) fails,
$$
(5.5) <(a^2, q_1) Q_1^{- 5/4} Q_1^{1/8}. R^2.\tag 5.11
$$
The contribution to (5.4) is at most
 $$
\frac{B^{2+\ve}}{PR} Q_1^3 Q_2^{5/2} \frac 1{|\Cal F|^2} \sum_{a\in\Cal A} |\Cal F_a|^2 \sum_{q\sim Q_1} (a^2, q) Q_1^{-\frac 98}R^2<
\frac {B^{2+\ve}R^3}{P} Q^{\frac {23} 8}< P^{\frac {15}{8}} R^7<\frac {P^2}{R^2}.\tag 5.12
$$

Finally, the contribution of $a=a_1$ and $(u, v, u_1, v_1)$ satisfying (5.9). We obtain
$$
\frac {B^{2+\ve}Q^2_1Q_2^{5/2}}{R.P} . \frac 1{|\Cal F|^2}  \Big(\frac P{Qd_1}\Big)^4 \sum_{q\sim Q_1} \ \sum_{a\in \Cal A} (a^2, q).\,(5.13)\tag 5.14
$$
with
$$
(5.13)=\Big|\Big\{(f, f_1, u, v, u_1, v_1)\in \Cal F_a^2\times \Big[1, \frac {Qd_1}P\Big]^4; f(u, v)=f(u_1, v_1)\Big\}\Big|.
$$
Recalling (1.3), it follows that
$$
(5.13)\ll R^\ve \frac {Q^2}{P^2} B^2|\Cal F_a|
$$
and
$$
\align
(5.14)&< R^\ve B^2 \frac {Q_1^2Q_2^{5/2}} {RP} \, \frac 1{|\Cal F|^2}  \, \frac {P^2}{Q^2} \sum_{q\sim Q_1} \ \sum_{a\in\Cal A} (a^2, q)|\Cal F_a|\\
&< B^2\frac {QP}{R^{1-\ve}|\Cal F|} \overset {(0.21)}\to < \frac {P^2}{R^{11/10}}\tag 5.15
\endalign
$$
(under proper assumption on $R$).

Hence from (5.2) and (5.15)
$$
(2.17) < \frac {P^2}R R^{-\frac 1{10}}.\tag 5.16
$$

\bigskip
\noindent
{\bf (6). Minor arcs estimate} 

From (4.15), (5.16), we obtain the following bound on (2.3)
$$
(2.3) <\frac {P^2}R \Big( R^{-\frac 1{10}}+Q^{-10^{-8}} \Big).\tag 6.1
$$
Take $R=(\log P)^{10^{10}}$ and $1\ll Q_0\leq R$.

Summing (6.1) over $Q>Q_0$ gives
$$
\int\limits_{\operatornamewithlimits\bigcup\limits_{q>Q_0} \operatornamewithlimits\bigcup\limits_{(b, q)=1} 
\Cal M(q, b)} \ | S_\omega (\theta)|^2 d\theta< \frac {P^2}{R} Q_0^{-10^{-8}}.\tag 6.3
$$
Returning to the definition of $\Cal M(q, b)$ in \S2,  we may further reduce the arcs $\Cal M(q, b)$, defining for $q<Q_0, (b, q)~=~1$
$$
\Cal M_0(q, b)=\Big[ \Big|\theta-\frac bq\Big|<\frac {Q_0^2}{RP^2}\Big].\tag 6.4
$$
It follows from (2.4) and (2.7) that for $\theta\in\Cal M(q, b), \theta=\frac bq+\vp$
$$
|S_\omega(\theta)| \lesssim \frac 1{\sqrt q R|\vp|} \tag 6.5
$$
and therefore
$$
\int\limits_{\Cal M (q, b)\backslash \Cal M_0(q, b)} |S_\omega (\theta)|^2 \lesssim \frac 1{qR^2} \int\limits_
{|\vp|>\frac{Q_0^2}{RP^2}} \frac 1{|\vp|^2} < \frac
{P^2}{RQ_0^2q}\tag 6.6
$$
which collected contribution is at most $\frac {P^2}{RQ_0}$.
Therefore
$$
\int\limits_{\Bbb T\backslash \operatornamewithlimits\bigcup\limits_{q\leq Q_0} \operatornamewithlimits\bigcup\limits_{(b, q)=1} 
\Cal M_0(q, b)} |S_\omega(\theta)|^2 d\theta<\frac {P^2}R \, Q_0^{-10^{-8}}.\tag 6.7
$$
where $1\ll Q_0 <(\log P)^{10^{10}}$ is a parameter.

\bigskip

\noindent
{\bf (7). Contribution of the major arcs}

Let
$$
1\ll Q_{00} <(\log P)^{\frac 14} \ \text { and }
Q_1 =\prod_{q|Q_{00}} q.\tag 7.1
$$
Let
$$
M=\frac {P^2}{Q_1}\tag 7.2
$$
and $\lambda$ a smooth even function on $\Bbb R$, supp\,$\lambda\subset [-1, 1], 0\leq\lambda \lesssim \frac 1M$ s.t.
$$
\sum_m \lambda\Big(\frac mM\Big)=1.\tag 7.3
$$
Let $\nu$ be the distribution on $\Bbb Z$ defined by
$$
\nu(n) =\sum_{m\in\Bbb Z}\lambda \Big(\frac mM\Big) \omega (n-m Q_1).\tag 7.4
$$
Hence
$$
\hat\nu(\theta) =\Big[\sum_m \lambda\Big(\frac mM\Big) e(-m Q_1\theta)\Big] S_\omega(\theta).\tag 7.4$'$
$$

Estimate
$$
\spreadlines{6pt}
\align
&\int_{\Bbb T} |\hat\nu (\theta) -S_\omega(\theta)|^2 d\theta =\int\Big| 1-\sum_m \lambda\Big(\frac mM\Big) e(-mQ_1\theta)\Big|^2
\ |S_\omega(\theta)|^2 d\theta\\
&< \frac {P^2}R Q_{00}^{-10^{-8}} +\sum_{q< Q_{00}} \, \sum_{(b, q)=1} \operatornamewithlimits\int\limits_{\Cal M_{00}(q, b)} 
\Big| 1-\sum\lambda
\Big(\frac mM\Big) e(-mQ_1\theta)\Big|^2 \ |S_\omega(\theta)|^2d\theta
\endalign
$$
where we applied (6.7) with $Q_0$ replaced by $Q_{00}$ and denote
$$
\Cal M_{00} (q, b) =\Big[\Big|\theta -\frac bq\Big|< \frac {Q^2_{00}}{RP^2}\Big].\tag 7.5
$$
Since $q|Q_1$ and $\lambda$ is even
$$
\spreadlines{6pt}
\align
&\operatornamewithlimits\int\limits_{\Cal M_{00} (q, b)}\Big| 1-\sum_m \lambda\Big(\frac mM\Big) e(-mQ_1\theta)\Big|^2 \ |S_\omega(\theta)|^2 d\theta\\
&\overset{(6.5)}\to \leq \operatornamewithlimits\int\limits_{|\vp|<\frac{Q_{00}^2}{RP^2}} \max_{|m|<M} |1-\cos (m Q_1\vp)|^2 \frac 1{q R^2\vp^2} d\vp\\
&<\frac {M^4Q_1^4}{qR^2} \Big(\frac{Q_{00}^2}{RP^2}\Big)^3 \overset {(7.2)}\to < \frac {P^2Q_{00}^6}{qR^5}
\endalign
$$
and we obtain from the choice of $R, Q_{00}$
$$
\Vert \hat\nu -S_\omega\Vert^2_2 < \frac {P^2}{R} Q_{00}^{-10^{-8}}.\tag 7.6
$$
Hence, by Parseval
$$
\sum |\nu (n) -\omega(n)|^2 <\frac {P^2}R Q_{00}^{-10^{-8}}.\tag 7.7
$$
Next, from definition of $\nu$ and $\omega$, clearly
$$
\nu(n) \sim\frac 1M \ \operatornamewithlimits{\overline\sum}_{\Cal F} |\{(x, y)\in [1, P]^2; (x, y)=1,  f_a(x, y) \equiv n(\mod Q_1) \text { and } \ |f_a(x, y)-n|<P^2\}|.
\tag 7.8
$$
Replacing the condition $(x, y)=1$ by the weaker condition $(x, y, p)=1$ for $p<B$ introduces an error at most $\frac {P^2}B$ with respect to the $\ell^1(\Bbb Z)$-norm.

With this replacement, we obtain the following lower bound on (7.8) (since \hbox{$Q_1, B\ll P^\ve)$}.
$$
\frac 1M \operatornamewithlimits{\overline\sum}\limits_{\Cal F} 
\Big|\Big\{ 0\leq x, y<Q_1; (x, y, Q_{00})=1 \text { and } f_a (x, y)\equiv n (\mod Q_1)\Big\}\Big|. \,(7.9)
$$
where (7.9) equals
$$
\align
&\text{mes\,} \Big[|s|, |t|<\frac P{Q_1}; \Big|f(s, t) -\frac n{Q_1^2}\Big| <\frac {P^2}{Q_1^2}\Big]\tag 7.10\\
&-\sum_{Q_{00} <p<B}\text{mes\,}\Big[|s|, |t|<\frac P{pQ_1}; \Big|f(s,t) -\frac n{p^2Q^2_1}\Big|<\frac{P^2}{p^2Q^2_1}\Big].\tag 7.10'
\endalign
$$
From the assumption on $f\in\Cal F$, it follows that taking $n\sim RP^2$
$$
\text{mes\,} \Big[ |s|, |t|<\frac P{Q_1}; \Big|f(s, t)-\frac n{Q_1^2}\Big|<\frac {P^2}{Q_1^2}\Big] \sim \frac {P^2}{Q_1^2 R}\tag 7.11
$$
and similarly for the terms in (7.10').
Hence $\nu(n)$ may be substituted by
$$
\spreadlines{6pt}
\align
\nu(n)&\sim \frac 1{Q_1R} \ \operatornamewithlimits{\overline\sum}_{f\in\Cal F}|\{0\leq x, y< Q_1; f(x, y)\equiv n+a (\mod Q_1)\ \text { and } (x, y, Q_{00})=1\}|\\
&=\frac 1R \ \operatornamewithlimits{\overline\sum}_{f\in\Cal F} \ \prod_{p<Q_{00}}\ [p^{-r} |\{0\leq x, y<p^r;  (x, y, p)=1 \ \text { and } \
f(x, y)\equiv n+a (\mod p^r)\}|]\tag 7.12
\endalign
$$
where $Q_1 =\prod_{p<Q_{00}} p^r$.

Recalling (0.24), $a\in\Cal A$ has no prime factors less than $R^{\frac 1{320}}$ and hence $(a, Q_1)=1$.

Let 
$$
f(x, y) =Ax^2+2Bxy+Cy^2 \in\Cal F \ \text { where } \  a^2=AC-B^2.
$$
Fixing $p< Q_{00}$, we may, since $(A, B, C)=1$, assume $(A, p)=1$
(the other cases are similar, except for $p=2$, $A, C$ even and $B$ odd; we leave the adjustment to the reader).

The equation
$$
f_a(x, y) \equiv n+a \ (\mod p^r)
$$
becomes
$$
(Ax+By)^2 +a^2y^2\equiv A(n+a) \ (\mod p^r)
$$
or equivalently
$$
x^2+y^2 \equiv A(n+a) \ (\mod p^r).\tag 7.13
$$
We seek for a lower bound on the number of solutions of (7.13) with $(x, y, p)=1$.

Assume $p>2$ (the local factor at $p=2$  requires an additional congruence assumption on $n$ at the place $p=2$).

Clearly, the number of solutions of (7.13) is at least
$$
p^{r-1}|\{(x, y)\in (\Bbb Z/ p\Bbb Z)^*\times \Bbb Z/p\Bbb Z; x^2+y^2 \equiv A(n+a) (\mod p)\}|.\tag 7.14
$$
Using Gauss sums, we obtain
$$
\spreadlines{6pt}
\align
&|\{(x, y)\in (\Bbb Z/p\Bbb Z)^*\times \Bbb Z/p\Bbb Z; x^2 + y^2 \equiv A(n+a) (\mod p)\}|\\
&=\frac 1p \ \sum^{p-1}_{b=0} \ \Big[\sum^{p-1}_{x=1} e_p (bx^2)\Big] \ \Big[ \sum^{p-1}_{y=0} e_p (by^2)\Big] e_p\big(-bA(n+a)\big)
\\
&= p-1 +\frac 1p \ \sum_{b=1}^{p-1} \Big(\zeta\sqrt p -\Big(\frac bp\Big)\Big) \zeta \sqrt p \, e_p\big(-bA(n+a)\big)\tag 7.15
\endalign
$$
with $\zeta =1$ (resp. $\zeta =i$) if $p\equiv 1(\mod 4)$ \big(resp. $p\equiv 3 (\mod 4)$\big).

Assume
$$
(n+a, p)= 1.\tag 7.16
$$
For $p\equiv 1(\mod 4)$, 
$$
(7.15) = p-2-\frac 1{\sqrt p} \, \sum^{p-1}_{b=1} \Big(\frac bp\Big) e_p\big(-bA(n+a)\big)= p-2+\alpha,
|\alpha| =1.
$$
For $p\equiv 3 (\mod 4)$,
$$
(7.15) = p-\frac i{\sqrt p} \, \sum_{b=1}^{p-1} \Big(\frac bp\Big) e_p \big(-bA(n+a)\big) =p+\alpha', |\alpha'|=1.
$$
Hence, if (7.16),
$$
(7.14)\geq \left\{\matrix p^r\Big(1-\frac 3p\Big) &\ \text { if \ $p\equiv 1 \ (\mod 4)$}\\
p^r\Big(1-\frac 1p\Big) & \ \text { if }  \ p\equiv 3 \ (\mod 4).\endmatrix \right. \tag 7.17
$$
Consequently
$$
(7.12) \gtrsim \frac 1R\ \operatornamewithlimits{\overline\sum}\limits_{\Cal F}{\!}{^{^*}} \ 1_{[(n+a, p)=1\text { for } 2< p<Q_{00}]} \ \prod_{3<p<Q_{00}} 
\Big(1-\frac 3p\Big)
 $$
and
$$
\nu(n) >(\log Q_{00})^{-C} R^{-1} \ \operatornamewithlimits{\overline \sum}_{\Cal F}{\!}^{^*} \ 1_{[(n+a, p) =1 \text { for } 2<p< Q_{00}]}\tag 7.18
$$
where $\sum^*$ refers to an additional congruence condition on $n, f$ $(\mod 8)$.

As an immediate consequence of (7.7) and (7.18), we obtain an alternative proof of the `positive density' conjecture, first established in [B-F1].

\proclaim
{Theorem 1}
Any integral Apollonian circle packing produces a set of curvatures of positive density in $\Bbb Z$.
\endproclaim

\noindent
{\bf Proof.}

With previous notation, let $S\subset [1, P^2 R]\cap \Bbb Z$ satisfy $|S|>(1-\tau)P^2R$ with $\tau >0$ a sufficiently small constant.
We show that
$$
\sum_{n\in S}\omega(n)>0.\tag 7.19
$$
Since $S$ is arbitrary, this will imply that $[1, P^2R]$ contains at least $\tau P^2R$ curvatures and hence Theorem 1.

Let $Q_{00}$ be a large constant and $K=\prod_{2<p<Q_{00}} p$.
From (7.7), (7.18)
$$
\align
\sum_{n\in S} \omega(n)&\geq \sum_{n\in S}\nu(n) -|S|^{1/2} \Big(\sum |\nu (n)-\omega(n)|^2\Big)^{1/2}\\
&>(\log Q_{00})^{-C} R^{-1} \operatornamewithlimits{\overline\sum}_{\Cal F} \underset {n\in S}\to {\sum\nolimits ^*} 1_{[(n+a, K)=1]} 
-P^2 Q_{00}^{-\frac 12 10^{-8}}\tag 7.20
\endalign
$$
(where $(*)$ refers to a congruence condition $(\mod 8)$ on $n+a)$.

Next
$$
\align
\underset{n\in S}\to{\sum\nolimits^*} 1_{[(n+a, K)=1]} &>\underset {n<P^2R}\to {\sum\nolimits^*} \ 1_{[(n+a, K)=1]} -\tau P^2 R\\
&\geq \frac 18 \prod_{2<p<Q_{00}} \ \Big(1-\frac 1p\Big) P^2R-\tau P^2 R\\
&\gtrsim (\log Q_{00})^{-1} P^2 R\tag 7.21
\endalign
$$
for $\tau=\tau (Q_{00})$ small enough.
Substituting in (7.20) gives
$$
\sum_{n\in S} \omega(n) >P^2 [c(\log Q_{00})^{-C} -Q_{00}^{-\frac 12 10^{-8}}]>0
$$
for $Q_{00}$ large enough.

\noindent
{\bf Remark.}
It should be noted that the previous argument establishes a stronger statement in fact.
It follows indeed that the set of curvatures in the ACP contains a coset of the integers, up to a `small' exceptional set (where `small' refers to small
density).
We do not attempt here to make quantitatively stronger statements, as better results will be obtained in the forthcoming paper [B-K2].

\bigskip

\noindent
{\bf (8). Prime Curvatures}

Another application of our analysis is an analogue of Theorem 1 for primes.  Thus

\proclaim
{Theorem 2}
The set of curvatures produced by any primitive integral Apollonian circle packing contains a subset of the primes of positive density.
\endproclaim

Denote  $\Cal P$ the set of primes and assume $S\subset [1, P^2R]\cap \Cal P$ satisfies
$$
\sum_{n\in S} \Lambda (n)>(1-\tau) P^2 R.\tag 8.1
$$
A straightforward adjustment of the proof of Theorem 1 shows that it will suffice to establish a bound
$$
\Big| \sum_{n\in S} \Lambda (n) [\omega(n) -\nu(n)]\Big|< Q_{00}^{-c} P^2\tag 8.2
$$
($c>0$ some fixed constant, $Q_{00}$ is a sufficiently large constant and $\tau=\tau(Q_{00})>0$\big).

Defining
$$
T(\theta) =\sum_{n\in S} \Lambda (n)e(n\theta)\tag 8.3
$$
we have to bound
$$
\Big|\int_{\Bbb T} [S_\omega -\hat\nu] T\Big | \leq \int_{\Bbb T} |S_\omega-\hat\nu|. |T|.\tag 8.4
$$
Invoking an additional ingredient, we will rely on certain distributional properties of exponential sums of the form (8.3), assuming $|S|\sim P^2R$.

As a consequence of the majorant property for the set of primes (see [B], [Gr]), there is the distributional inequality
$$
\text{mes\,}[\theta\in\Bbb T; |T(\theta)|>\delta P^2 R]\ll_\ve \delta^{-2-\ve} (P^2 R)^{-1} \text { for } \delta> (\log P)^{-A}
\tag 8.5
$$
(and hence for all $\delta>0$).

Hence, if $\Omega\subset\Bbb T$ such that $P^2R|\Omega| >1$, we have
$$
\int_\Omega |T(\theta)|^2 d\theta\ll (|\Omega|P^2 R)^\ve P^2 R.\tag 8.6
$$
 
Take $Q_0 =(\log P)^{10^{10}}$ and let $\Cal M_0(q, b), \Cal M_{00}(q, b)$ be defined by (6.4), (7.5).

Decompose (8.4) as follows
$$
\spreadlines{6pt}
\align
|(8.4)|&\leq\int_{\Bbb T\backslash \operatornamewithlimits\bigcup\limits_{\Sb q\leq Q_0\\ (b, q) =1\endSb} \Cal M_0 (q, b)}
 \ |S_\omega -\hat\nu | \ |T|\tag 8.7\\
& +\sum_{\Sb Q_{00} <Q<Q_0\\ Q\text{ dyadic}\endSb} \ \int\limits_{\operatornamewithlimits\bigcup\limits_{\Sb q\sim Q\\ (b, q)=1\endSb} 
\Cal M_0(q, b)} \ |S_\omega -\hat\nu| \ |T|\tag 8.8\\
&+\sum_{\Sb q\leq Q_{00}\\ (b, q)=1\endSb} \int_{\Cal M_0(q, b)} \ |S_\omega -\hat\nu| \, |T|\tag 8.9
\endalign
$$
Recall that $|\hat\nu|\leq |S_\omega|$ by (7.4$'$).
\medskip

Since certainly $\Vert T\Vert_2 <(\log P)^{\frac 12} \sqrt R P$, we have
$$
(8.7) \lesssim (\log P)^{\frac 12} \sqrt R P\Biggl[ \int\limits_{\Bbb T\backslash \operatornamewithlimits\bigcup\limits_{q\lesssim Q_0} 
\Cal M_0(q, b)} |S_\omega|^2\Biggr]
^{1/2} \overset{(6.7)}\to < (\log P)^{\frac 12} \sqrt R.P.\frac P{\sqrt R} \, Q_0^{-\frac 12 10^{-7}} <P^2(\log P)^{-4}.\tag 8.10
$$
We estimate (8.8).
Fix $Q_{00} <Q<Q_0$ ($Q$ dyadic) and define
$$
\Cal M'(q, b)=\Bigl[\Big|\theta -\frac bq \Big|< \frac {q^2}{RP^2}\Big].\tag 8.11
$$
Then
$$
\spreadlines{6pt}
\align
&\operatornamewithlimits\int_{\operatornamewithlimits\bigcup\limits_{\Sb q\sim Q\\ (b, q)=1\endSb}\Cal M_0(q, b)} |S_\omega-\hat\nu| \ |T|\leq\\
&\operatornamewithlimits\int_{\operatornamewithlimits\bigcup\limits_{\Sb q\sim Q\\ (b, q)=1\endSb}\Cal M'(q, b)} 
|S_\omega| . |T|+\sum_{\Sb q\sim Q\\ (b, q)=1\endSb} \  
\operatornamewithlimits\int\limits_{\Cal M_0 (q. b)\backslash \Cal M' (q, b)} |S_\omega|.|T| = (8.12) + (8.13).
\endalign
$$
By (6.7) and (8.6)
$$
(8.12) < \frac P{\sqrt R} Q^{-\frac 12 10^{-7}}\Vert T|_{_{\operatornamewithlimits\bigcup\limits_{q\sim Q}\Cal M' (q, b)}}\Vert_2 
\ll \frac P{\sqrt R} \, Q^{-\frac 12 10^{-7}}
Q^{3\ve} P\sqrt R < Q^{-\frac 13 10^{-7}} P^2\tag 8.14
$$
and, using (6.5) and (8.6)
$$
\spreadlines{6pt}
\align
(8.13) &\lesssim \sum_{\Sb q\sim Q\\ (b, q)=1\endSb} \frac 1{RQ^{1/2}} \operatornamewithlimits\int\limits_{\frac {Q^2}{RP^2}<|\vp|<\frac {Q_0^2}{RP^2}}
|\vp|^{-1} |T\Big(\frac bq +\vp\Big)\Big|{d\vp}\\
&< \frac 1{RQ^{1/2}} \ \sum_{s, Q^2< 2^s< Q^2_0} P^2R{2^{-s}} \sum_{\Sb q\sim Q\\(b, q) =1\endSb} \ \operatornamewithlimits\int\limits_{|\vp|\sim\frac {2^s}{RP^2}}
\Big| T\Big(\frac bq +\vp\Big)\Big|d\vp\\
&\ll \frac {P^2}{Q^{1/2}} \sum_{2^s>Q^2} 2^{-s} \Big( \frac {2^s}{RP^2}Q^2\Big)^{\frac 12} (2^s Q^2)^\ve P\sqrt R< P^2 Q^{-1/3}.\tag 8.15
\endalign
$$
and from (8.14), (8.15)
$$
(8.8) < Q_{00}^{-\frac 13 10^{-7}} P^2.\tag 8.16
$$

Finally,
$$
\spreadlines{6pt}
\align
(8.9) &\leq \operatornamewithlimits\int\limits_{\operatornamewithlimits\bigcup\limits_{q\leq Q_{00}}
\ \operatornamewithlimits\bigcup\limits_{(b, q)=1} [|\theta-\frac bq|< \frac {Q_{00}^2}{RP^2}]}
|S_\omega - \hat\nu| \, |T| +\sum_{\Sb q\leq Q_{00}\\ (b, q)=1\endSb} \
 \operatornamewithlimits\int\limits_{ [\frac {Q^2_{00}}{RP^2} \leq |\theta -\frac bq|<\frac {Q_0^2}{RP^2}]}
|S_\omega| \, |T|\\
&= (8.17)+(8.18)
\endalign
$$
From (7.6), (8.6),
$$
(8.17) \ll \Vert S_\omega -\hat\nu\Vert_2 \ Q_{00}^{3\ve} \, (P^2 R)^{1/2} <\frac P{\sqrt R} \, Q_{00}^{-\frac 12 10^{-7}} \, 
Q_{00}^{3\ve} \, P\sqrt R < Q_{00}^{-\frac 13
10^{-7}} P^2\tag 8.19
$$
and (8.18) is bounded as (8.13)
$$
(8.18) <P^2 Q_{00}^{-1/3}.\tag 8.20
$$
Hence
$$
(8.9) < Q_{00}^{-\frac 13 10^{-7}}P^2.\tag 8.21
$$
From (8.10), (8.16), (8.21),
$$
(8.4)< Q_{00}^{-\frac 13 10^{-7}} P^2.\tag 8.22
$$
which is the desired inequality (8.2).
This proves Theorem 2.


\Refs
\widestnumber\no{XXXXXXXX}

\ref\no{[B]} \by J.~Bourgain
\paper On $\Lambda(p)$-subsets of the squares
\jour Israel J.~Math. 67 (1989), 291--311
\endref


\ref\no{[BEW]}\by B.~Berndt, R.~Evans, K.~Williams
\paper Gauss and Jacobi sums
\jour Canadian Math. Soc., Vol. 21 (1998)
\endref

\ref\no{[B-F1]}\by J.~Bourgain, E.~Fuchs
\paper A proof of the positive density conjecture for integer Apollonian circle packings
\jour to appear in JAMS
\endref

\ref\no{[B-F2]}\by J.~Bourgain, E.~Fuchs
\paper On representation of integers by binary quadratic forms
\jour submitted to IMRN
\endref


\ref\no{[BGS]}\by J.~Bourgain, A.~Gamburd, P.~Sarnak
\paper Generalization of Selberg's 3/16 theorem and affine sieve
\jour to appear in Acta Math
\endref

\ref\no{[B-K]}\by J.~Bourgain, A.~Kontorovich
\paper
On representations of integers in thin subgroups of $SL_2(\Bbb Z)$
\jour GAFA, Vol. 20, N5, 2010, p. 1144--1174
\endref

\ref\no{[B-K2]} \by J.~Bourgain, A.~Kontorovich
\paper On the strong density conjecture for integral Apollonian circle packings
\jour in preparation
\endref



\ref\no{[F-S]}\by E.~Fuchs, K.~Sanden
\paper
Prime numbers and the local to global principle in Apollonian circle packings
\jour preprint
\endref


\ref\no{[G]}\by B.~Green
\paper Roth's theorem in the primes
\jour Annals of Math. (2), 161 (2005), no 3, 1609--1636
\endref

\ref\no{[G-S]}\by V.~Guillemin, S.~Sternberg
\paper Geometric asymptotics
\jour Math.~Surveys 14, AMS 1977
\endref


\ref\no{[GLMWY]}\by R.~Graham, L.~Lagarias, C.~Mallows, A.~Wilks, C.~Yan
\paper Apollonian circle packings: number theory
\jour J. of Number Theory, 100 (1), pp 1--45 (2003)
\endref

\ref\no{[K-O]}\by A.~Kontorovich, H.~Oh
\paper Apollonian circle packings and closed horospheres on hyperbolic 3-manifolds
\jour to appear in JAMS
\endref

\ref\no{[S]} \by P.~Sarnak
\paper Letter to Lagarias on integral Apollonian packings
\jour (2007)
\endref

\endRefs
\enddocument